\documentclass[12pt,reqno,a4wide]{amsart}

\usepackage{tikz} 

\usetikzlibrary{decorations.pathreplacing}
\usetikzlibrary{arrows,decorations.pathmorphing,snakes}
\usetikzlibrary{fit}
\usepackage{calrsfs}

\allowdisplaybreaks

\begin{document}
	
\title {Generating functions for a lattice path model introduced by Deutsch}

\author[Helmut Prodinger]{Helmut Prodinger}
\address{Department of Mathematics, University of Stellenbosch 7602,
	Stellenbosch, South Africa}
\email{hproding@sun.ac.za}

\thanks{During the Corona crisis, I don't have access to a library. A more polished version of the current material will hopefully be posted soon.}

\begin{abstract} The lattice path model suggested by E. Deutsch is derived from ordinary Dyck paths, but
	with additional down-steps of size $-3,-5,-7,\dots$\;. For such paths, we find the generating functions of
	them, according to length, ending at level $i$, both, when considering them from left to right and from right
	to left. The generating functions are intrinsically cubic, and thus (for $i=0$) in bijection to various objects, like
	even trees, ternary trees, etc.
	
\end{abstract}

\maketitle

\section{Introduction}
Nonnegative lattice paths consisting of up-steps $(1,1)$ and down-steps $(1,-1)$ and ending at the $x$-axis again are enumerated by Catalan numbers. Richard Stanley \cite{Stanley} devoted a whole bookon combinatorial objects enumerated by Catalan numbers.

Nonnegative lattice paths consisting of up-steps $(1,1)$ and down-steps $(1,-2)$ and ending at the $x$-axis again are enumerated by 
the numbers $\frac1{2n+1}\binom{3n}{n}$. These generalized Catalan numbers appear less often than Catalan numbers but they are sequence $A001764$ in \cite{OEIS} and also enumerate many different combinatorial objects. It is recommended to watch Donald Knuth's christmas lecture \cite{christ}. 

Deutsch~\cite{Deutsch} suggested different down-steps: $(1,-1)$, $(1,-3)$, $(1,-5)$, \dots. This still leads to the numbers
$\frac1{2n+1}\binom{3n}{n}$, and, because of the plethora of combinatorial objects, allows for bijective methods.
Indeed, Cameron~\cite{Cameron00} established a bijection to so-called even trees, which are plane trees with the property that each node has
an \emph{even} number of subtrees.

In this paper, we will construct generating functions for the enumeration of ``Deutsch-paths''. We also generalize the concept by
allowing that such a path ends at level $i\ge0$, $i=0$ being the original model.

Since up-steps and down-steps are no longer symmetric (as in ordinary Dyck paths), the enumeration from right to left, working with
up-steps $(1,1)$, $(1,3)$, $(1,5)$, \dots\ and one down-step $(1,-1)$ leads (for $i\neq0$) for different generating functions and thus enumerations. This right-to-left model is more difficult to treat.

The method consists by first assuming a horizontal barrier at level $h$, so that the Deutsch-paths live in a strip. Then the generating functions are rational, and linear algebra methods (Cramer's rule) can be employed. In the resulting formul\ae\ it is then possible to perform the limit $h\to\infty$.

To know the enumeration both, from left to right and from right to left allows to compute the total area, summed over all Deutsch-paths returning to the $x$-axis. The area of a path is the sum of all its ordinates.

Similar methods were already used in \cite{Prodinger-Cameron}, solving an open problem in Cameron's thesis~\cite{Cameron}. However, the results from our earlier paper cannot be used for the present enumeration of Deutsch-paths.

In whatever follows, the variable $z$ is used in generating functions to mark the length of the paths. The substitution
$z^2=t(1-t)^2$ is used throughout, always with this meaning. In this way, the cubic equations that we encounter become manageable. Compare this with \cite{ProSelWag}.

We find explicit expressions for the generating functions $f_i(z)$ and $g_i(z)$ of Deutsch-paths ending on level $i$ from left to right resp.\ from right to left. They are rational in an auxiliary variable $t$ which in itself satisfies a cubic equation. Thus one might speak about cubic Deutsch-paths. We are also able to write down explicit expressions for the coefficients of these generating functions.
As a corollary we find an explicit expression for the cumulative area, summed over all Deutsch-paths of a given length. This result is derived from the generating functions. It should be stated that computer algebra (Maple) played a major rule in the current project.

\section{Two generating functions}

The following expansions will be used later. Set
\begin{align*}
\frac{1}{1-X+z^2X^3}=\sum_{n\ge0}a_nX^n
\end{align*}
and
\begin{align*}
r_1&=1-t,\qquad r_{2,3}=\frac{t\pm\sqrt{4t-3t^2}}{2}.
\end{align*}
We sometimes find it useful to abbreviate $W=\sqrt{4t-3t^2}$.
A direct computation confirms that
\begin{equation*}
(1-r_1X)(1-r_2X)(1-r_3X)=1-X+z^2X^3.
\end{equation*}
Some background information why this cubic equations factors nicely can be found in \cite{christ, ProSelWag}.
Then, adjusting the initial values, we get the explicit expression
\begin{equation*}
a_n=\frac{1}{3t-1}\Bigl[-r_1^{n+1}+\frac{3t+W}{2W}r_2^{n+1}- \frac{3t-W}{2W}r_3^{n+1}\Bigr].
\end{equation*}

The other expansion we need is related to
\begin{align*}
\frac{1}{1-Y^2-zY^3}=\sum_{n\ge0}b_nY^n.
\end{align*}
It can be checked directly via $(1-\mu_1 z)(1-\mu_2 z)(1-\mu_3 z)=1-Y^2-zY^3$ that
\begin{equation*}
b_n=\frac{1}{3t-1}[A\mu_1^n+B\mu_2^n+C\mu_3^n]
\end{equation*}
with
\begin{align*}
\mu_1=\frac{z}{t-1}, \qquad\mu_2=\frac{-z(t+W)}{2t(t-1)},\qquad \mu_3=\frac{z(-t+W)}{2t(t-1)}
\end{align*}
and
\begin{equation*}
A=t,\qquad B=\frac{2t-1}{2}+\frac{t}{2W},\qquad C=\frac{2t-1}{2}-\frac{t}{2W}.
\end{equation*}
The quantities  $A$, $B$, $C$ take care of the initial values.

\section{Enumeration of   paths from left to right}

\begin{center}
	\begin{tikzpicture}[scale=0.3]
	
	\draw[step=1.cm,black] (-0.0,-0.0) grid (16,4.0);

	\draw[ultra thick] (0,0) to (16,0);

	\draw[thick] (0,0) -- (1,1) -- (2,2)-- (3,1)-- (4,2)-- (5,3)-- (6,2)-- (7,3)-- (8,2)--(9,1) --(10,2)--(11,3)--(12,0)
	--(13,1)--(14,2)--(15,3)--(16,2);
	
	\node at (-1.5,0.0){$(0,0)$};
	
	\end{tikzpicture}
	\end	{center}

The picture shows a Deutsch  path ending in $(16,2)$ and being bounded by $4$ (or higher).

Let $a_{n,k}$ be the number of Deutsch-paths ending at $(n,k)$ and being bounded by $h$. 
In order not to clutter the notation, we did not put the letter $h$ into the definition, especially, 
since $h$ has only a very temporary meaning here. 

The recursion (for $n\ge1$)
\begin{equation*}
a_{n,k}=a_{n-1,k-1}+a_{n-1,k+1}+a_{n-1,k+3}+a_{n-1,k+5}+\cdots
\end{equation*}
with the understanding that $a_{n,k}$ should be interpreted as 0 if $k<0$ or $k>h$ is easy to understand. The starting value
is $a_{0,0}=1$. It is natural to introduce the generating functions
\begin{equation*}
f_k=f_k(z)=\sum_{n\ge0}a_{n,k}z^n.
\end{equation*}
Here is a little list:
\begin{align*}
f_0&=1+{z}^{2}+3\,{z}^{4}+12\,{z}^{6}+55\,{z}^{8}+268\,{z}^{10}+1338\,{z}^{
	12}+6741\,{z}^{14}+34075\,{z}^{16}+\cdots\\
f_1&=z+2\,{z}^{3}+7\,{z}^{5}+30\,{z}^{7}+142\,{z}^{9}+701\,{z}^{11}+3517\,{
	z}^{13}+17751\,{z}^{15}+\cdots\\
f_2&={z}^{2}+3\,{z}^{4}+12\,{z}^{6}+55\,{z}^{8}+268\,{z}^{10}+1338\,{z}^{12
}+6741\,{z}^{14}+34075\,{z}^{16}+\cdots\\
f_3&={z}^{3}+4\,{z}^{5}+18\,{z}^{7}+87\,{z}^{9}+433\,{z}^{11}+2179\,{z}^{13
}+11010\,{z}^{15}+\cdots\\
f_4&={z}^{4}+5\,{z}^{6}+25\,{z}^{8}+126\,{z}^{10}+637\,{z}^{12}+3224\,{z}^{
	14}+16324\,{z}^{16}+\cdots\\
f_5&={z}^{5}+6\,{z}^{7}+32\,{z}^{9}+165\,{z}^{11}+841\,{z}^{13}+4269\,{z}^{
	15}+\cdots\\
f_6&={z}^{6}+7\,{z}^{8}+39\,{z}^{10}+204\,{z}^{12}+1045\,{z}^{14}+5314\,{z}
^{16}+\cdots
\end{align*}
The recursion for the numbers $a_{n,k}$ translates into
$$f_k=zf_{k-1}+zf_{k+1}+zf_{k+3}+zf_{k+5}+\cdots+[\![k=0]\!],$$ which is best written as a matrix equation
\begin{equation*}
\begin{bmatrix}
1& -z&0&-z&0&-z&0&\dots\\
-z&1&-z&0 &-z&0&-z&\dots\\
0&-z&1&-z&0&-z&0&\dots\\
&&&\vdots\\
&&&&&& -z&1
\end{bmatrix}
\begin{bmatrix}
f_0\\f_1\\f_2\\ \vdots\\f_h
\end{bmatrix}
=\begin{bmatrix}
1\\0\\0\\ \vdots\\0
\end{bmatrix}
\end{equation*}

Now let $d_m$ be the determinant of this matrix with $m$ rows and columns. We have $d_0=1$, $d_1=1$, $d_2=1-z^2$, and the recursion
\begin{equation*}
d_m=d_{m-1}-z^2d_{m-3}.
\end{equation*}
The characteristic equation of this recursion is the cubic equation
\begin{equation*}
\lambda^3-\lambda^2+z^2=0.
\end{equation*}
Note also the generating function
\begin{align*}
R(X)&=\sum_{m\ge0}d_{m-1}X^m=1+X+X^2+\sum_{m\ge3}(d_{m-2}-z^2d_{m-4})X^m\\
&=1+X+X^2+X\sum_{m\ge2}d_{m-1}X^m-z^2X^3R(X)\\
&=1+X+X^2+XR(X) -X-X^2-z^2X^3R(X),
\end{align*}
or
\begin{equation*}
R(X)=\frac{1}{1-X+z^2X^3}=\sum_{j\ge0}d_{j-1}X^j.
\end{equation*}
So we see that $d_{n-1}$ are the numbers studied in the previous section.
Therefore
\begin{equation*}
d_m=\frac{1}{3t-1}\Bigl[-r_1^{m+2}+\frac{3t+W}{2W}r_2^{m+2}- \frac{3t-W}{2W}r_3^{m+2}\Bigr].
\end{equation*}
Cramer's rule now leads to
\begin{equation*}
f_k= z^{k}\frac{d_{h-k}}{d_{h+1}},
\end{equation*}
which, when performing the limit $h\to\infty$, leads to
\begin{equation*}
f_k(z)= z^{k}r_1^{-k-1}=\frac{z^{k}}{(1-t)^{k+1}}.
\end{equation*}
This form will be useful later, but we would also like to compute $[z^n]f_k(z)$, i.~e., the numbers $a_{n,k}$. We can only have contributions (which is also clear for combinatorial reasons) if $n\equiv k \bmod 2$. So let us set $n=2N+i$, $k=2K+i$ for 
$i=0,1$, and compute
\begin{align*}
[z^{2N+i}]\frac{ z^{2K+i}}{(1-t)^{2K+i+1}}&=[z^{2N-2K}]\frac{1}{(1-t)^{2K+i+1}}.
\end{align*}
It helps, as mentioned before, to set $z^2=x=t(1-t)^2$. Then we can continue
\begin{align*}
a_{2N+i,2K+i}&=[x^{N-K}]\frac{1}{(1-t)^{2K+i+1}}\\
&=\frac1{2\pi i}\oint\frac{dx}{x^{N-K+1}}\frac{1}{(1-t)^{2K+i+1}}\\
&=\frac1{2\pi i}\oint\frac{dt(1-t)(1-3t)}{t^{N-K+1}(1-t)^{2N-2K+2}}\frac{1}{(1-t)^{2K+i+1}}\\
&=[t^{N-K}](1-3t)\frac{1}{(1-t)^{2N+i+2}}\\
&=\binom{3N-K+i+1}{N-K}-3\binom{3N-K+i}{N-K-1}.
\end{align*}
Notice in particular the enumeration of paths ending at the $x$-axis:
\begin{align*}
a_{2N,0}
&=\binom{3N+1}{N}-3\binom{3N}{N-1}=\frac1{2N+1}\binom{3N}{N},
\end{align*}
a generalized Catalan number, enumerating many different combinatorial objects, as mentioned in the Introduction.

\section{Enumeration of ternary paths from right to left}

We still prefer to work from left to right, so we change our setting as follows:
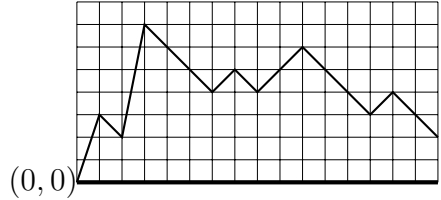
\begin{center}
	\begin{tikzpicture}[scale=0.3]
	
	\draw[step=1.cm,black,ultra thin] (-0.0,-0.0) grid (16.0,8.0);

	\draw[ultra thick] (0,0) to (16,0);

	\draw[thick] (0,0) -- (1,3) -- (2,2)-- (3,7)-- (4,6)-- (5,5)-- (6,4)-- (7,5)-- (8,4)
	-- (9,5)-- (10,6)-- (11,5)-- (12,4)-- (13,3)-- (14,4)-- (15,3)-- (16,2) ;
	
	\node at (-1.5,0.0){$(0,0)$};
	
	\end{tikzpicture}
\end	{center}
	The picture shows a reversed Deutsch  path ending in $(16,2)$ and being bounded by $7$ (or higher).

For the notation, we switch from $a_{n,k}$ to $b_{n,k}$ and from $f_k(z)$ to $g_k(z)$, and give a short list as an illustration.
\begin{align*}
g_0&=1+{z}^{2}+3\,{z}^{4}+12\,{z}^{6}+55\,{z}^{8}+268\,{z}^{10}+1338\,{z}^{
	12}+6741\,{z}^{14}+34075\,{z}^{16}
+\cdots\\
g_1&=z+3\,{z}^{3}+12\,{z}^{5}+55\,{z}^{7}+273\,{z}^{9}+1428\,{z}^{11}+7752
\,{z}^{13}+43263\,{z}^{15}+\cdots\\
g_2&=2\,{z}^{2}+9\,{z}^{4}+43\,{z}^{6}+218\,{z}^{8}+1155\,{z}^{10}+6324\,{z
}^{12}+35511\,{z}^{14}
+\cdots\\
g_3&=z+6\,{z}^{3}+31\,{z}^{5}+163\,{z}^{7}+882\,{z}^{9}+48967\,{z}^{11}+
27759\,{z}^{13}+\cdots\\
g_4&=3\,{z}^{2}+19\,{z}^{4}+108\,{z}^{6}+609\,{z}^{8}+3468\,{z}^{10}+20007
\,{z}^{12}
+\cdots\\
g_5&=z+10\,{z}^{3}+65\,{z}^{5}+391\,{z}^{7}+2313\,{z}^{9}+13683\,{z}^{11}+\cdots\\
g_6&=4\,{z}^{2}+34\,{z}^{4}+228\,{z}^{6}+1431\,{z}^{8}+8787\,{z}^{10}
+\cdots\\
\end{align*}

The linear system changes now as follows:

\begin{equation*}
\begin{bmatrix}
1& -z&\dots\\
-z&1& -z&\dots\\
0&-z&1& -z&\dots\\
-z&0&-z&1& -z&\dots\\
0&-z&0&-z&1& -z&\dots\\
&&&\vdots\\
&&&&& -z&1
\end{bmatrix}
\begin{bmatrix}
g_0\\g_1\\g_2\\g_3\\g_4\\ \vdots\\g_h
\end{bmatrix}
=\begin{bmatrix}
1\\0\\0\\0\\0\\ \vdots\\0
\end{bmatrix}
\end{equation*}

The determinant of the matrix is the same as before by transposition: $d_{h+1}$. However, the application of Cramer's rule is more involved now. We must evaluate the determinant of
\begin{align*}
&\begin{bmatrix}
0& 0&\dots& 1 &\dots&0\\
-z&1& -z&0&\dots\\
0&-z&1& 0&\dots\\
-z&0&-z&0& 0&\dots\\
&&&\vdots\\
&&&0& 0&1
\end{bmatrix}\\
&\ \underbrace{\phantom{xxxxxxxxxxxx}}\quad\underbrace{\phantom{xxxxxx}}\\
&\hspace*{1.2cm}q-1\hspace*{1.6cm}m-q
\end{align*}
  Call this determinant $\Delta_{m,q}$. We want to find a recursion for it.

By expansion, we find the recursions
\begin{equation*}
\alpha_{m}=\alpha_{m-1}-z^2\alpha_{m-3}
\end{equation*}
for $\alpha_{m}=\Delta_{m,q}$, for fixed $q$, and
\begin{equation*}
\beta_{m}=\beta_{m-2}+z\beta_{m-3}
\end{equation*}
for $\beta_{m}=\Delta_{m,m-q}$, for fixed $q$.

Eventually, with a lot of help by \textsf{Gfun}, we find that for $2\le q<m$:
\begin{align*}
\Delta_{m,q}&=[X^mY^q]
\frac{zX^2Y^2(1+zXY+zYX^2+z^2Y^2X^3)}{(1-X+z^2X^3)(1-X^2Y^2-zX^3Y^3)}\\
&=z[X^{m-q}Y^{q-2}]\frac{(1+zY)(1+zXY)}{(1-X+z^2X^3)(1-Y^2-zY^3)},\\
\Delta_{m,1}&=d_{m-1}=[X^{m-1}]\frac{1-z^2X^2}{1-X+z^2X^3},\\
\Delta_{m,m}&=z^2[Y^{m-2}]\frac{1+zY}{1-Y^2-zY^3}.
\end{align*}
We can now continue with the computation, according to Cramer's rule:
\begin{align*}
g_i&=\frac{\Delta_{h+1,i+1}}{d_{h+1}}.
\end{align*}

The instance $\Delta_{m,m}$ is less interesting, since we will eventually push the barrier $h$ to infinity, and
then it plays no role anymore. The instance $\Delta_{m,1}$ will show as that $f_0=g_0$, which is clear from combinatorial
reasons, since it describes the same objects when left and right are switched. So we concentrate on the remaining cases
and compute some  generating functions:
\begin{align*}
\Delta_{m,q}
&=z[X^{m-q}Y^{q-2}]\frac{(1+zY)(1+zXY)}{(1-X+z^2X^3)(1-Y^2-zY^3)}\\*
&=za_{m-q}[Y^{q-2}]\frac{1+zY}{1-Y^2-zY^3}
+z^2a_{m-q-1}[Y^{q-3}]\frac{1+zY}{1-Y^2-zY^3}\\*
&=za_{m-q}(b_{q-2}+zb_{q-3})+z^2a_{m-q-1}(b_{q-3}+zb_{q-4})
\end{align*}
and we will rewrite this according to $g_i=\frac{\Delta_{h+1,i+1}}{d_{h+1}}$.

At that stage, we will let $h\to\infty$ and assume that $i\ge1$.
\begin{align*}
g_i&=\frac{z}{(1-t)^{i+2}}(b_{i-1}+zb_{i-2})+\frac{z^2}{(1-t)^{i+3}}(b_{i-2}+zb_{i-3})\\
&=\frac1{3t-1}A\mu_1^{i-3}\Big[\frac{z\mu_1^2}{(1-t)^{i+2}}+\frac{z^2\mu_1}{(1-t)^{i+2}}
+\frac{z^2\mu_1}{(1-t)^{i+3}}+\frac{z^3}{(1-t)^{i+3}}\Big]\\
&+\frac1{3t-1}B\mu_2^{i-3}\Big[\frac{z\mu_2^2}{(1-t)^{i+2}}+\frac{z^2\mu_2}{(1-t)^{i+2}}
+\frac{z^2\mu_2}{(1-t)^{i+3}}+\frac{z^3}{(1-t)^{i+3}}\Big]\\
&+\frac1{3t-1}C\mu_3^{i-3}\Big[\frac{z\mu_3^2}{(1-t)^{i+2}}+\frac{z^2\mu_3}{(1-t)^{i+2}}
+\frac{z^2\mu_3}{(1-t)^{i+3}}+\frac{z^3}{(1-t)^{i+3}}\Big]\\
&=\frac{t+(2t-1)W}{2(3t-1)(1-t)^{i+1}W}\mu_2^{i-3}\Big[\frac{z\mu_2^2}{1-t}+\frac{z^2\mu_2}{1-t}
+\frac{z^2\mu_2}{(1-t)^2}+\frac{z^3}{(1-t)^2}\Big]\\
&+\frac{-t+(2t-1)W}{2(3t-1)(1-t)^{i+1}W}\mu_3^{i-3}\Big[\frac{z\mu_3^2}{1-t}+\frac{z^2\mu_3}{1-t}
+\frac{z^2\mu_3}{(1-t)^2}+\frac{z^3}{(1-t)^2}\Big]\\
&=z\frac{t+(2t-1)W}{4(3t-1)(1-t)^{i+2}W}\mu_2^{i-3}(-3t^2+3t+2-(t-3)W)\\
&+z\frac{-t+(2t-1)W}{4(3t-1)(1-t)^{i+2}W}\mu_3^{i-3}(-3t^2+3t+2+(t-3)W)\\
&=\frac{z}{2(1-t)^{i+2}W}\mu_2^{i-3}(t^3-5t^2+5t-(t^2-t-1)W)\\
&+\frac{z}{2(1-t)^{i+2}W}\mu_3^{i-3}(-t^3+5t^2-5t-(t^2-t-1)W)\\
&=\frac{1}{4W(1-t)^{i+5}}\mu_2^{i}(t^3-5t^2+5t-(t^2-t-1)W)(2t^2-3t+W)\\
&+\frac{1}{4W(1-t)^{i+5}}\mu_3^{i}(-t^3+5t^2-5t-(t^2-t-1)W)(2t^2-3t-W)\\
&=\frac{t}{2(1-t)^{i+1}W}\mu_2^{i}(-t+2+W)+\frac{t}{2(1-t)^{i+1}W}\mu_3^{i}(t-2+W)\\
&=\frac{t}{2(1-t)^{i+1}}(\mu_2^{i}+\mu_3^{i})-\frac{z(t-2)}{2(1-t)^{i+2}}\frac{\mu_2^{i}-\mu_3^{i}}{\mu_2-\mu_3}.
\end{align*}
So we managed to compute the functions $g_i$ for all $i\ge0$.

\section{Extraction of coefficients}

We start by noticing the pleasant formulae
\begin{equation*}
\mu_2\mu_3=t-1,\quad \mu_2+\mu_3=\frac{z}{1-t}.
\end{equation*}
We need  the Girard-Waring formula, see e. g. \cite{Gould}:
\begin{equation*}
X^m+Y^m=\sum_{0\le k\le m/2}(-1)^k\binom{m-k}{k}\frac{m}{m-k}(XY)^{k}(X+Y)^{m-2k}.
\end{equation*}
In our instance, we need
\begin{equation*}
\mu_2^i+\mu_3^i=\sum_{0\le k\le i/2}z^{i-2k}\binom{i-k}{k}\frac{i}{i-k}(1-t)^{3k-i} .
\end{equation*}
The other version is
\begin{equation*}
\frac{X^m-Y^m}{X-Y}=\sum_{0\le k\le (m-1)/2}(-1)^k \binom{m-1-k}{k}(XY)^{k}(X+Y)^{m-1-2k}
\end{equation*}
and in the special instance
\begin{equation*}
\frac{\mu_2^i-\mu_3^i}{\mu_2-\mu_3}=\sum_{0\le k\le (i-1)/2}z^{i-1-2k}\binom{i-1-k}{k}(1-t)^{3k+1-i}.
\end{equation*}
So we get for $i\ge1$
\begin{align*}
g_i&=\frac{t}{2(1-t)^{i+1}}(\mu_2^{i}+\mu_3^{i})-\frac{z(t-2)}{2(1-t)^{i+2}}\frac{\mu_2^{i}-\mu_3^{i}}{\mu_2-\mu_3}\\
&=\frac{t}{2}\sum_{0\le k\le i/2}z^{i-2k}\binom{i-k}{k}\frac{i}{i-k}(1-t)^{3k-2i+1}\\
&-\frac{t-2}{2}\sum_{0\le k\le (i-1)/2}z^{i-2k}\binom{i-1-k}{k}(1-t)^{3k-2i+1}\\
&=t\sum_{1\le k\le (i-1)/2}z^{i-2k}\binom{i-1-k}{k-1}(1-t)^{3k-2i+1}\\
&+\sum_{0\le k\le (i-1)/2}z^{i-2k}\binom{i-1-k}{k}(1-t)^{3k-2i+1}.
\end{align*}
Now we distinguish the two cases $i$ even resp.\ odd. Set $i=2I+\delta$, with $\delta\in\{0,1\}$.
\begin{align*}
	[z^{2N+\delta}]g_i
	&=[z^{2N}]t\sum_{1\le k\le (i-1)/2}z^{2I-2k}\binom{i-1-k}{k-1}(1-t)^{3k-2i+1}\\*
	&+[z^{2N}]\sum_{0\le k\le (i-1)/2}z^{2I-2k}\binom{i-1-k}{k}(1-t)^{3k-2i+1}\\
	&=[x^{N}]\sum_{1\le k\le (i-1)/2}t^{I-k+1}\binom{i-1-k}{k-1}(1-t)^{k-i+1-\delta}\\*
	&+[x^{N}]\sum_{0\le k\le (i-1)/2}t^{I-k}\binom{i-1-k}{k}(1-t)^{k-i+1-\delta}\\
	&=\sum_{1\le k\le (i-1)/2}[t^{N}](1-t)^{-2N-1}(1-3t)t^{I-k+1}\binom{i-1-k}{k-1}(1-t)^{k-i+1-\delta}\\*
	&+\sum_{0\le k\le (i-1)/2}[t^{N}](1-t)^{-2N-1}(1-3t)t^{I-k}\binom{i-1-k}{k}(1-t)^{k-i+1-\delta}\\
	&=\sum_{1\le k\le (i-1)/2}\binom{i-1-k}{k-1}[t^{N-I+k-1}](1-3t)(1-t)^{-n+k-i}\\*
	&+\sum_{0\le k\le (i-1)/2}\binom{i-1-k}{k}[t^{N-I+k}](1-3t)(1-t)^{-n+k-i}\\
	&=\sum_{1\le k\le (i-1)/2}\binom{i-1-k}{k-1}\biggl[\binom{3N+2\delta+I-2}{N-I+k-1}-3\binom{3N+2\delta+I-3}{N-I+k-2}\biggr]\\*
	&+\sum_{0\le k\le (i-1)/2}\binom{i-1-k}{k}\biggl[\binom{3N+2\delta+I-1}{N-I+k}-3\binom{3N+2\delta+I-2}{N-I+k-1}\biggr]\\
\end{align*}

\section{The area  } 

Each contribution $c_i$ to the area of a path $(0,c_0=0),\dots,(2n,c_{2n}=0)$ can be seen as splitting the  path into a path of
length $k$ (left to right) ending at level $i$ and a path of length $2n-k$ (right to left) also ending at level $i$.
Since we are working with generating functions, all possible such splittings are taken into account when taking the product of two such generating functions.

The cumulated area is thus given as (write again $z^2=x=t(1-t)^2$)
\begin{align*}
\textsc{area}&=\sum_{i\ge0}if_i(z)g_i(z)=\frac{t(1+3t)}{(1-t)(1-3t)^2}.
\end{align*} 
Then
\begin{align*}
[z^{2n}]\textsc{area}&=[x^n]\frac{t(1+3t)}{(1-t)(1-3t)^2}\\
&=\frac1{2\pi i}\oint \frac{dx}{x^{n+1}}\frac{t(1+3t)}{(1-t)(1-3t)^2}\\
&=\frac1{2\pi i}\oint \frac{dt(1-t)(1-3t) }{t^{n+1}(1-t)^{2n+2}}\frac{t(1+3t)}{(1-t)(1-3t)^2}\\
&=[t^{n-1}]\frac{1}{(1-t)^{2n+2}}\frac{1+3t}{1-3t}\\
&=\sum_{k\ge0}3^k[t^{n-1-k}]\frac{1+3t}{(1-t)^{2n+2}}\\
&=\sum_{k\ge0}3^k\bigg[\binom{3n-k}{n-1-k}+3\binom{3n-1-k}{n-2-k}\biggr].
\end{align*} 

The paper \cite{BG} contains general results about the area of lattice paths, of
a less explicit nature than what we are doing here. Anyway, since Deutsch-paths 
have an infinite set of possible steps, they do not fall into the framework studied
in \cite{BG}.

\bigskip

We hope that the gentle reader will find our analysis of cubic Deutsch-paths exciting;
we are sure that there are many questions left for future research, possibly also 
of a bijective nature.

\bibliographystyle{plain}


\end{document}